\documentclass[12pt]{article}
\makeindex
\usepackage[hypertex]{hyperref}
\usepackage{amsfonts}
\usepackage{amsmath}
\usepackage[all]{xy}
\usepackage{mathrsfs}

\newtheorem{myle}{{\sc Lemma}}[section]

\newtheorem{mypr}{{\sc Proposition}}[section]

\newtheorem{myre}{{\sc Remark}}[section]

\newcommand{\mysection}[1]{\section{#1}\setcounter{equation}
{0}}
\begin{document}
\textwidth 5in
\textheight 7in
\begin{center}
{\Large\bf Corrigendum to "The $q$-analogue of bosons and Hall algebras" and some remarks}
\end{center}
\begin{center}
Youjun Tan
\\
\vspace{2mm}
{\small Mathematical college, Sichuan University, Chengdu, 610064, P.R.China\\
ytan@scu.edu.cn}
\end{center}
MSC 2000: 16W35, 17B37
\mysection{Acknowledgements and backgrounds.}
 I thank Akira Masuoka very much for the following reasons. At first I am very grateful to him for his comments in \cite{1} that the proof of \cite[Proposition 3.2 (3)]{2} is false: The argument that $0=g(u_{1})=Xg(u_{2})=X$ in \cite[line -3, p.4346]{2} is wrong, since $X$ needs not to be in ${\rm\bf B}^{+}(\Lambda)$. In fact, the statement \cite[Proposition 3.2 (3)]{2} itself is wrong: There indeed exists an object in ${\cal O}({\rm\bf B})$ which is {\it not} semisimple if the $q$-boson ${\rm\bf B}$ is determined by a Borcherds-Cartan (or generalized Kac-Moody) form on an {\it infinite} set (a counter example is given below), although there is only one isoclass of simple objects in ${\cal O}({\rm\bf B})$.  

Secondly, although I proved that the main statement, \cite[Theorem 3.1]{2}, due to B. Sevenhant and M. Van den Bergh, can be deduced from \cite[Proposition 3.2 (1)]{2} and \cite[Proposition 3.1]{2}, as Masuoka pointed out to me, \cite[Theorem 3.1]{2} can be proved much more directly
 using his result in \cite{1}.

Thirdly, as Masuoka pointed out to me, even the proof of \cite[Proposition 3.2 (1)]{2} can be simplified to large extent by using the natural skew-pairing on  ${\rm\bf U}^{-} \otimes {\rm\bf U}^{+}$ described in \cite{1}. Moreover, he also pointed out to me that the argument in \cite[page 4345]{2} may lead confusion between the left multiplication by  ${\rm F}_i\in {\rm\bf U}^{+}$ and the natural action by ${\rm F}_i$ on ${\rm\bf B}^{+}(\Lambda)$ ($= {\rm\bf U}^{+}$). In deed, all formula in \cite[page 4345]{2} such that ${\rm F}_{i}P_{2}=0={\rm F}_{i}Q_{2}$, ${\rm F}_{i}X_{j}=0$, etc., should mean that ${\rm F}_{i}P_{2}=P'_{2}{\rm F}_{i}$, ${\rm F}_{i}Q_{2}=Q_{2}'{\rm F}_{i}$, ${\rm F}_{i}X_{j}=X_{j}'{\rm F}_{i}$ in ${\rm\bf U}$, etc., where $P'_{2}$, $Q_{2}'$ and $X_{j}'$ belong to ${\rm\bf U}^{+}$. For details please see {\sc Remark} \ref{Lemma3.4} below.

In the last month I've been trying to seek a "correct" proof of \cite[Proposition 3.2 (3)]{2}. It is Masuoka who always finds mistakes in those arguments. I feel sorry for wasting so much time of him.

\hspace{1mm}\newline The counter example given below is motived by investigation of the semisimplicity of ${\cal O}({\rm\bf B})$ for the case of $q$-boson ${\rm\bf B}$ determined by a Borcherds-Cartan form on a {\it finite} set. Note that in this case the semisimplicity follows by using extremal projectors: The Kac-Moody case is due to Nakashima, while the more general case is due to Masuoka, see {\sc Remark} \ref{Extremalprojector} below. Moreover, Masuoka's generalized extremal projectors deduces a nontrivial semisimple subcategory of ${\cal O}({\rm\bf B})$ in the case of infinite indexed set, for details see \cite[Theorem 4.4]{1}.

\hspace{1mm}\newline For self-contained purpose, we keep the the following notations. Let ${\cal I}$ be a countable set. A Borcherds-Cartan form on
${\cal I}$ is a non-degenerate ${\rm\bf Q}$-valued bilinear form (-,-)
satisfying the following conditions (a)-(c):

(a) (-,-) is symmetric;

(b) $(i,j)\leq 0$ for $i,j\in {\cal I}$ if $i\not=j$ and

(c) $\frac{2(i,j)}{(i,i)}$ is an integer if $(i,i)$ is positive.

The elements of ${\cal I}$ are called simple roots and we have a disjoint union
${\cal I}={\cal I}^{{\rm re}}\cup {\cal I}^{{\rm im}}$ where
${\cal I}^{{\rm re}}$ (resp.  ${\cal I}^{{\rm im}}$) contains the
elements $i\in {\cal I}$ such that $(i,i)$ $>0$ (resp. $(i,i)\leq 0$).
For a real root $i,$ we set $a_{ij}=-2\frac{(i,j)}{(i,i)},$ and $d_{i}=\frac{(i,i)}{2}$,  $q_{i}=q^{d_{i}}$, where $q$ is fixed to be an indeterminant.

\hspace{1mm}\newline By definition, the $q$-boson, also called Kashiwara algebra, ${\rm\bf B}$  associated to the Borcherds-Cartan form on ${\cal I}$ is
an associative algebra over ${\rm\bf Q}(q)$ generated by symbols
${\rm E}_{i},{\rm F}_{i}$ for $i\in{\cal I}$ subject to the following
relations (\ref{def1})-(\ref{def4}):
\begin{eqnarray}\label{def1}
{\rm F}_{i}{\rm E}_{j}=q^{(i,j)}{\rm E}_{j}{\rm F}_{i}+\delta_{ij}\,\,\,
{\rm for}\,\,\, i,\,j\,\,{\rm in}\,\, {\cal I};
\end{eqnarray}
\begin{eqnarray}\label{def2}
\sum_{t=0}^{a_{ij}+1}(-1)^{t}\left [\begin{array}{c}a_{ij}+1\\t\end{array}
\right ]_{d_{i}}{\rm E}_{i}^{t}{\rm E}_{j}{\rm E}_{i}^{a_{ij}+1-t}=0\,\,{\rm for\,\, real\,\, simple\,\, root}\,\, i,
\end{eqnarray}
\begin{eqnarray}\label{def3}
\sum_{t=0}^{a_{ij}+1}(-1)^{t}\left [\begin{array}{c}a_{ij}+1\\t\end{array}
\right ]_{d_{i}}{\rm F}_{i}^{t}{\rm F}_{j}{\rm F}_{i}^{a_{ij}+1-t}=0\,\,{\rm for\,\, real\,\, simple\,\, root}\,\, i,
\end{eqnarray}
\begin{eqnarray}\label{def4}
{\rm E}_{i}{\rm E}_{j}-{\rm E}_{j}{\rm E}_{i}=0,\,\,\,{\rm F}_{i}{\rm F}_{j}
-{\rm F}_{j}{\rm F}_{i}=0\,\,{\rm for\,\, any\,\, pair}\,\, i,j\,\,{\rm with}\,\, (i,j)=0,
\end{eqnarray}
where $[\,]_{d_{i}}$ is the standard notation of quantum binomials.
\begin{myre}  For a symmetrizable Kac-Moody algebra $\mathfrak{g}$, the
$q$-boson ${\rm\bf B}_{q}(\mathfrak{g})$ is defined in {\rm \cite[3.3]{0}}. Here we adopt the ``positive'' version.
\end{myre}
Following Kashiwara \cite{0}, we define ${\cal O}({\rm\bf B})$ to be the category containing left ${\rm\bf B}$-modules $M$ such that for any element $u$ of $M$
there exists an integer $l$ with ${\rm F}_{i_{1}}{\rm F}_{i_{2}}
\ldots{\rm F}_{i_{l}}u=0$ for any $i_{1},$ $i_{2},\ldots,i_{l}$ in ${\cal I}$. Note that the category ${\cal O}({\rm\bf B})$ is closed under subs, quotients and extensions. Thus ${\cal O}({\rm\bf B})$ is an abelian subcategory of the category of left ${\rm\bf B}$-modules.

Let ${\rm\bf B}^{+}$ (resp. ${\rm\bf B}^{-}$) be the subalgebra of ${\rm\bf B}$ generated by ${\rm E}_{i}$ (reps. ${\rm F}_{i}$), $i\in{\cal I}$. Since ${\rm\bf B}={\rm\bf B}^{+}{\rm\bf B}^{-}$, due to (\ref{def1}), the {\it Verma} module ${\rm\bf B}/{\rm\bf B}^{-}$ is isomorphic to ${\rm\bf B}^{+}$ with module structure given by ${\rm F}_{i}1=0$ for all $i\in {\cal I}$. Then we have the following
\begin{myle}\label{SimpleOb} {\rm (\cite[Proposition 3.2 (1)]{2}).}
${\rm\bf B}^{+}$ is a simple object of ${\cal O}({\rm\bf B})$. Moreover, ${\rm\bf B}^{+}$ represents the unique isoclass of simple objects in ${\cal O}({\rm\bf B})$. \hfill $\Box$
\end{myle}
\begin{myre}\label{Lemma3.4}
As mentioned above, this result can be obtained by using Masuoka's result in {\rm \cite{1}}. My proof depends on {\rm \cite[Lemma 3.4]{2}}, and the formula
\begin{eqnarray*}
{\rm F}_{i}^{a}{\rm E}_{i}^{a}Z=\frac{1-q^{(a-1)(i,i)}}{1-q^{(i,i)}}{\rm F}_{i}^{a-1}{\rm E}_{i}^{a-1}Z=\ldots=\prod_{t=1}^{a-1}\frac{1-q^{t(i,i)}}{1-q^{(i,i)}}Z
\end{eqnarray*}
(Note also that $(i,i)=0$ may appear) in {\rm \cite[page 4345]{2}} should be replaced by the formula in ${\rm\bf U}$:
$${\rm F}_{i}^{a}{\rm E}_{i}^{a}Z$$
\begin{eqnarray}\label{Z}
=q^{(a-1)(i,i)}{\rm F}_{i}^{a-1}{\rm E}_{i}^{a-1}Z_{i}'{\rm F}_{i}+(1+q^{(i,i)}+\ldots +q^{(a-2)(i,i)}){\rm F}_{i}^{a-1}{\rm E}_{i}^{a-1}Z,
\end{eqnarray}
where ${\rm F}_{i}Z=Z_{i}'{\rm F}_{i}$ for some $Z_{i}'\in{\rm\bf U}^{+}$. Since ${\rm F}_{i}X_{1}=X'_{1}{\rm F}_{i}$ for some $X_{1}'\in {\rm\bf U}^{+}$ and ${\rm F}_{i}u_{\lambda}=0$, applying the action of ${\rm F}_{i}^{l_{i}}$ to ${\rm E}_{i}^{l_{i}}X_{1}u_{\lambda}$ it follows that
\begin{eqnarray*}
{\rm F}_{i}^{l_{i}}{\rm E}_{i}^{l_{i}}X_{1}u_{\lambda}=(1+q^{(i,i)}+\ldots +q^{(l_{i}-2)(i,i)}){\rm F}_{i}^{\l_{i}-1}{\rm E}_{i}^{l_{i}-1}X_{1}u_{\lambda},
\end{eqnarray*}
and, if $b>a$ then ${\rm F}_{i}^{b}{\rm E}_{i}^{a}Xu_{\lambda}=0$, whenever ${\rm F}_{i}X=X'{\rm F}_{i}$. (For a more general expression see {\rm \cite[(6.4)]{-1}}\,). The remaining argument goes through and {\rm \cite[Lemma 3.4]{2}} follows. \hfill $\Box$
\end{myre}
We have the following 
\begin{myle}\label{cyclicsimple}
If a nonzero cyclic module ${\rm\bf B}m\in {\cal O}({\rm\bf B})$ is simple then ${\rm\bf B}m={\rm\bf B}^{+}m$ as vector spaces. If a nonzero cyclic module ${\rm\bf B}m$ satisfies that ${\rm F}_{i}m=0$ for all $i\in {\cal I}$, then ${\rm\bf B}m$ is simple and ${\rm\bf B}m={\rm\bf B}^{+}m$.
\end{myle}
{\bf Proof.} Assume that ${\rm\bf B}m$ is simple. Since ${\rm\bf B}m\in{\cal O}({\rm\bf B})$, there is a $Y\in{\rm\bf B}^{-}$ such that $Ym\not=0$ but ${\rm F}_{j}Ym=0$ for all $j\in {\cal I}$. By {\sc Lemma} \ref{SimpleOb} it follows that ${\rm\bf B}m={\rm\bf B}^{+}Ym$, which is $\mathbb{Z}_{+}{\cal I}$-graded. Thus there is a unique $X\in{\rm\bf B}^{+}$ such that $m=XYm$. If $Pm=Qm$ for some $P,Q\in{\rm\bf B}^{+}$, then $PXYm=QXYm$, which means that $PX=QX\in{\rm\bf B}^{+}$ by \cite[Lemma 3.4]{2}. Therefore $P=Q$ and hence there is an isomorphism of vector spaces ${\rm\bf B}^{+}m\simeq {\rm\bf B}^{+}Ym$ induced by $m\mapsto Ym$. So ${\rm\bf B}^{+}m={\rm\bf B}m$ as required. \hfill $\Box$

Note that ${\rm\bf B}^{+}$ is $\mathbb{Z}_{+}{\cal I}$-graded as a ${\rm\bf Q}(q)$-module:
\begin{eqnarray}\label{B+}
{\rm\bf B}^{+}=\oplus_{\alpha\in\mathbb{Z}_{+}{\cal I}}{\rm\bf B}^{+}_{\alpha}\,,
\end{eqnarray}
where ${\rm\bf B}^{+}_{\alpha}$ is spanned by the monomials of the form ${\rm E}_{i_{1}}\ldots {\rm E}_{i_{t}}$ with $i_{1}+\ldots +i_{t}=\alpha$.
\mysection{Cases of finite indexed sets.}
Let us recall the following
\begin{mypr}\label{Propofinite}
{\rm (Kashiwara-Masuoka-Nakashima).} Assume that the indexed set ${\cal I}$ is finite. Then the category ${\cal O}({\rm\bf B})$ is semisimple, that is, every nonzero object in ${\cal O}({\rm\bf B})$ is a sum of simple objects, and hence isomorphic to a sum of copies of ${\rm\bf B}^{+}$. \hfill $\Box$
\end{mypr}
\begin{myre}\label{Extremalprojector}  For the $q$-boson ${\rm\bf B}_{q}(\mathfrak{g})$ associated to a symmtrizable Kac-Moody algebra $\mathfrak{g}$, M. Kashiwara stated firstly that ${\cal O}({\rm\bf B}_{q}(\mathfrak{g}))$ is semisimple in {\rm \cite{0}} without explicit proof.  T. Nakashima {\rm \cite{1.5}} proved that there is a well defined element $\Gamma$ in some completion of ${\rm\bf B}_{q}(\mathfrak{g})$, called the {\rm extremal projector}, satisfying that
\begin{eqnarray}\label{extremal}
\begin{split}
{\rm F}_{i}\Gamma=\Gamma{\rm E}_{i}=0,\,\,\Gamma^{2}=\Gamma,\\
\sum\limits_{k\geq 0} a_{k}\Gamma b_{k}=1\,\,\,{\rm for\,\,some}\,\,a_{k}\in {\rm\bf B}_{q}^{+}(\mathfrak{g}),\,\,b_{k}\in{\rm\bf B}_{q}^{-}(\mathfrak{g}).
\end{split}
\end{eqnarray}
Applying the action of $\Gamma$, Nakashima proved the semisimplicity of ${\cal O}({\rm\bf B}_{q}(\mathfrak{g}))$. Masuoka generalized this construction to a more general situation, including the case of $q$-boson associated to symmetrizable Borcherds-Cartan form {\rm \cite[Proposition 3.6]{1}}. These constructions generalize the rank 1 case due to Kashiwara {\rm \cite{0}} (see also {\rm \cite{-1}}) in a remarkable and highly nontrivial way.\hfill $\Box$
\end{myre}
\begin{myre}
Assume that ${\cal I}$ is infinite. In {\rm \cite{1}} Masuoka considered a subcategory ${\cal O}'({\rm\bf B})$ of left ${\rm\bf B}$-modules $M$ such that
\begin{description}
  \item {\rm (1)} $M$ is an object of ${\cal O}({\rm\bf B})$.
  \item {\rm (2)} For any $m\in M$, there is a finite set $F(m)$ such that ${\rm F}_{i_{1}}\ldots{\rm F}_{i_{t}}\ldots {\rm F}_{i_{r}}m=0$ for any $i_{t}\not\in F(m)$.
\end{description}
(See {\rm \cite[Definition 4.2]{1}}). Notations in {\rm \cite{1}} is adjusted here for brevity. Then the subcategory ${\cal O}'({\rm\bf B})$ is shown by Masuoka to be equivalent to $Vec$, which means that it is semisimple. Thus in this case Masuoka's generalized extremal projector is crucial in my view.
\end{myre}
\mysection{A counter example to the case of infinite indexed sets.} Assume that ${\cal I}=\{0,1,2,\ldots\}$ is infinite. For any sequence $\{a_{j}\}_{j\geq 1}$ with $0\not=a_{j}\in {\rm\bf Q}(q)$, set
\begin{eqnarray}\label{N}
N={\rm\bf B}/J,\,\,J\,\,{\rm is\,\,the\,\,left\,\,ideal\,\,generated\,\,by}\,\,{\rm F}_{j}-a_{j}{\rm F}_{0}:\,j\geq 1,\,{\rm F}_{0}^{2}.
\end{eqnarray}
Then $N$ becomes a left ${\rm\bf B}$-module in a natural way. Clearly $N$ is a nonzero object of ${\cal O}({\rm\bf B})$. Let $u\in N$ be the image of $1\in{\rm\bf B}$. By definition in $N$ it holds that
\begin{eqnarray}\label{idenN}
{\rm F}_{j}u=a_{j}{\rm F}_{0}u;\,\,\,{\rm F}_{r}{\rm F}_{s}u=0,\,\,j\geq 1,r,s\in{\cal I}.
\end{eqnarray}
Note that $N$ has a decomposition as vector spaces:
\begin{eqnarray}\label{Ndecompo}
N={\rm\bf B}^{+}u\oplus {\rm\bf B}^{+}{\rm F}_{0}u,
\end{eqnarray}
where ${\rm\bf B}^{+}{\rm F}_{0}u$ is simple by {\sc Lemma} \ref{cyclicsimple}, since ${\rm F}_{j}{\rm F}_{0}u=0$ for all $j\in {\cal I}$.

We claim that $N$ is {\it not} semisimple. Assume contrarily that $N$ is semisimple. Then, by $\mathbb{Z}_{+}{\cal I}$-gradation there is a short exact sequence of the form
\begin{eqnarray}
\xymatrix{0\ar[r]&{\rm\bf B}^{+}{\rm F}_{0}u\ar[r]^-{f}&N\ar[r]^-{g}&{\rm\bf B}^{+}\ar[r]&0}, \label{sequence1}
\end{eqnarray}
which must split. It follows that $N$ has a simple submodule of the form ${\rm\bf B}^{+}(u+Q{\rm F}_{0}u)$ for some $Q\in {\rm\bf B}^{+}$. By (\ref{def1}), for all $j\geq 1$ it holds that in ${\rm\bf B}$:
\begin{eqnarray}\label{Qj}
{\rm F}_{j}Q=Q_{j}{\rm F}_{j}+Q_{j}^{'}:\,\,\,Q_{j},Q_{j}^{'}\in{\rm\bf B}^{+}.
\end{eqnarray}
Thus, for any $j\geq 1$, by (\ref{idenN}) and (\ref{Qj}) it follows that
\begin{eqnarray*}
0&=&{\rm F}_{j}(u+Q{\rm F}_{0}u)={\rm F}_{j}u+Q_{j}{\rm F}_{j}{\rm F}_{0}u+Q_{j}^{'}{\rm F}_{0}u\\
&=&{\rm F}_{j}u+Q_{j}^{'}{\rm F}_{0}u=(a_{j}+Q_{j}^{'}){\rm F}_{0}u,
\end{eqnarray*}
which means that $0\not=Q_{j}^{'}=-a_{j}\in {\rm\bf Q}(q)$ for all $j\geq 1$. But this is impossible in ${\rm\bf B}$, since ${\cal I}$ is infinite, there is always a $t\geq 1$ such that ${\rm E}_{t}$ does not appear in $Q$, and hence ${\rm F}_{t}Q=f_{t}(q)Q{\rm F}_{t}$ for some $f_{t}(q)\in{\rm\bf Q}(q)$ by (\ref{def1}), which implies that $a_{t}=0$, a contradiction. Therefore $N$ is not semisimple as claimed.

\end{document}